\documentclass[sn-mathphys-num]{sn-jnl}

\usepackage{graphicx}%
\usepackage{multirow}%
\usepackage{amsmath,amssymb,amsfonts}%
\usepackage{amsthm}%
\usepackage{mathrsfs}%
\usepackage[title]{appendix}%
\usepackage{xcolor}%
\usepackage{textcomp}%
\usepackage{manyfoot}%
\usepackage{booktabs}%
\usepackage{algorithm}%
\usepackage{algorithmicx}%
\usepackage{algpseudocode}%
\usepackage{listings}%
\usepackage{amssymb}
\usepackage{amsthm}
\usepackage{setspace}
\usepackage{amsmath}
\usepackage{graphicx}
\usepackage{MnSymbol}
\usepackage{cutwin}
\usepackage{lipsum}
\usepackage{soul,color}
\usepackage{subfigure}
\usepackage{caption}
\usepackage{subcaption}

\usepackage{geometry}
\begin{document}

 \title[Modeling pathway signal propagation]{Mathematical modeling of biochemical signal propagation in many-stage enzymatic pathways}

\author*{\fnm{Chathranee} \sur{Jayathilaka}}\email{chathraneeaththanapolaarachchlage@monash.edu}

\author{\fnm{Mark} B. \sur{Flegg}}\email{mark.flegg@monash.edu}
\equalcont{These authors contributed equally to this work.}


\affil{\orgdiv{School of Mathematics}, \orgname{Monash University}, \orgaddress{ \city{Clayton}, \postcode{3168}, \state{VIC}, \country{Australia}}}



\abstract{Biochemical signalling cascades transduce extracellular stimuli into cellular responses through sequences of discrete, node-to-node activations. While signal fidelity depends critically on local interaction kinetics, the mechanisms governing information propagation in realistic, highly variable kinetic contexts remain poorly understood. In this paper, we develop a mathematical framework for travelling waves in canonical feed-forward pathways governed by nonlinear Michaelis-Menten-type kinetics. For uniform pathways, we characterise the complete steady-state landscape and demonstrate that activation bias (the contribution of the binary states of each node to downstream activation) between connected nodes acts as a key bifurcation parameter dictating wave existence. Extending this framework to heterogeneous networks, we show how parameter gradients and random kinetic variations distort wavefronts and induce heavy fluctuations in propagation speed. To recover predictable signal transmission, we introduce a novel reciprocal-velocity spatial rescaling technique. We demonstrate that this coordinate transformation inherently absorbs local kinetic variations, effectively smoothing wave velocities and preserving wavefront profiles without requiring bespoke parameter tuning or continuous limits. Finally, by testing the framework's limits against extreme parameter variability, we reveal how severe kinetic bottlenecks lead to functional pathway fragmentation, offering a mathematically justified basis for rational model reduction in complex biochemical networks.}




\maketitle

\section{Introduction}
\label{sec:Intro}

Signal transduction is the fundamental process by which cells convert environmental stimuli into specific functional responses, a mechanism essential for sensing, survival, growth, and adaptation. This process typically involves complex networks of chemical interactions where primary signals, such as extracellular ligands binding to membrane-bound receptors, initiate a cascade of downstream molecular modifications. These biochemical cascades ultimately translate external information into precise cellular actions, including gene expression, cell cycle regulation, and differentiation \cite{hunter2000signaling}.

Biochemical signaling pathways frequently exhibit intricate dynamical behaviors, such as temporal oscillations and multistability \cite{Othmer1997SignalTransduction, Schwen2021NonlinearSignaling}. A particularly prevalent feature in biological transduction is the emergence of traveling waves---sequential pulses of activity that propagate through a medium or along a signaling cascade. Notable biological examples include intracellular calcium waves, the propagation of action potentials in neural networks, and the cAMP waves that coordinate the aggregation of \textit{Dictyostelium} \cite{martiel1987model}.

A central driver of these diverse dynamical behaviors is the intrinsic delay associated with signal propagation. In pathways characterized by negative feedback, for instance, sufficient intrinsic delay can precipitate oscillatory behavior. While multistep dynamical models consisting of linear transitions are well-known to converge to robust delays in the limit \cite{hurtado2019generalizations}, the delay properties of linear cascades exhibiting more complex chemical kinetics, such as those investigated in this study, have received comparatively less attention. This research gap is largely due to the lack of analytical tractability and the inherent variability in signal transduction mechanisms.

The functional necessity of these delays is underscored in numerous biological systems where specific behaviors are mathematically dependent on the finite rate of cascade propagation. Within mathematical frameworks, biochemical delays may be explicitly represented using delay-differential equations (DDEs), as seen in various gene regulation models \cite{Glass2021NonlinearDelayMotifs}. Alternatively, delays are often defined intrinsically through increased signaling complexity or the inclusion of intermediary steps, representing a mechanistically modeled delay \cite{Murray2002MathematicalBiologyI, McKeithan1995KineticProofreading}. Fixed-delay models generally imply that a downstream component is influenced after a long sequence of processes, effectively a chain of ordinary differential equations (ODEs); thus, the choice between explicit DDEs or implicit ODE chains is frequently a matter of modeling strategy, numerical efficiency or tractability.

Signal propagation through biological pathways is mediated by various mechanisms, including phosphorylation cascades, protein-protein interactions, and the diffusion of second messengers. Phosphorylation cascades, such as the MAPK/ERK pathway, exemplify how extracellular stimuli trigger sequential kinase activations that ultimately regulate nuclear gene expression \cite{pearson2001mitogen}. Similarly, the Wnt pathway transmits signals through protein complexes that modulate transcriptional activity, where the spatial and temporal characteristics of the pathway dictate how signals are interpreted \cite{clevers2012wnt}. ODE models serve as a standard quantitative tool for describing these processes, offering a structured approach to analyzing how molecular concentrations evolve over time in response to diverse stimuli.

A critical application of ODE models lies in the stability analysis of signal propagation. Such analysis characterizes how systems respond to perturbations and identifies attracting states, revealing whether a signal will decay, oscillate, or persist at a steady state. For example, Huang and Ferrell \cite{huang1996ultrasensitivity} demonstrated that ODE frameworks can effectively capture the ultrasensitive, switch-like responses of the MAPK/ERK pathway, enabling cells to transition decisively between signaling states. In systems biology, evaluating the stability of 'no-input' models is vital for understanding how cells maintain a baseline readiness while preventing premature activation. This stability analysis of resting states has significant implications for understanding the prevention of aberrant pathway activation in pathologies such as cancer and neurodegenerative diseases \cite{kholodenko2000negative}. Furthermore, by leveraging the relationship between DDEs and systems of ODEs, the stability properties of the former can often be more readily investigated through the latter \cite{scarabel2024equations}.

While linear (mass-action) transition chains provide a useful abstraction for studying delays, real signaling pathways are composed of heterogeneous biochemical interactions rather than identical, repeated steps. The effect of realistic variations in nonlinear enzyme kinetics across sequential stages of the chain can be obscured if pathways are treated as uniform chains. Rios \textit{et al.} \cite{Rios2021} used a generalised framework of enzymatic interactions governed by Michaelis--Menten kinetics to illustrate how network structure and kinetic nonlinearity shape responses like adaptation. However, as the complexity of biochemical networks increases, they quickly become analytically intractable, and reduced representations that preserve core dynamical behaviours are often heuristically accepted. In previous work, we examined and justified such reductions by introducing the \textit{oriented form} of a biochemical pathway---a formal diagrammatic representation which heuristically normalises the direction and sign of interactions whilst keeping the kinetic parameters fixed \cite{jayathilaka2024two}. That study demonstrated that networks/pathways sharing the same oriented form behave strictly the same only under highly symmetric and biologically restrictive assumptions but in under particular conditions the oriented pathways are good approximations for the adopted pathway it heuristically represents, and in others it does a very poor job. It is surprising how in many cases complex pathways can be simplified using the oriented form despite nonlinear chemical kinetics.

These findings motivate the present study, which examines these useful oriented pathways to gain insight into signal propagation speed and, by extension, the delay times associated with chains of biochemical cascades. Understanding the velocity and reliability of signal propagation is central to coordinating cellular responses. In multistep pathways, transmission delays can generate coherent activity fronts that advance through cascades, justifying our use of traveling-wave models describing transitions from activation and inactivation states. While signal propagation in continuous systems has been rigorously characterized by reaction--diffusion equations, such as the Fisher--KPP equation \cite{fisher1937wave, kolmogorov1982study}, intracellular signaling typically proceeds through discrete reaction steps. In these discrete systems, propagation speed is determined by the interplay of nonlinear kinetics, saturation phenomena, and spatial localization, as evidenced in MAPK cascades \cite{huang1996ultrasensitivity} and engineered signaling circuits \cite{o2011tunable}. Whilst we do not consider physical space in this manuscript, we aim to investigate the propagation speed within such biochemical cascades to characterize the temporal dynamics of signal transmission.

In this study, we apply a discrete traveling-wave framework to feed-forward `oriented' chemical pathways governed by nonlinear kinetics. We first analyze idealized cascades with identical, repeated chemical steps to establish a baseline for propagation. We then generalize this framework to account for pathway heterogeneity, where each node is characterized by its propensity for `on' and `off' states, and each edge $i$ is defined by distinct kinetic parameters, $\alpha_i$ and $\beta_i$. To quantify how this variability shapes the overall transmission delay, we introduce a velocity-based rescaling technique similar to reciprocal-velocity rescaling (RVR) \cite{aronson2006nonlinear}. This approach allows us to treat wave speed as an explicit parameter by defining distance metrics between nodes that are determined by local kinetic parameters. By transforming the governing equations into a coordinate frame that moves with the signal front, we isolate the intrinsic propagation properties of the oriented form from external scaling effects. This analytical strategy is complemented by sensitivity analysis, which identifies the specific reaction rates that most strongly govern wave dynamics, offering deeper mechanistic insights into how minor parameter fluctuations, such as those observed in the Ras--MAPK pathway \cite{heinrich2002mathematical} can induce substantial changes in signal velocity.

The structure of this paper is organized as follows. Section 2 establishes the mathematical framework for modeling signal propagation, extending the oriented form representation to general node-to-node interactions. We begin by analyzing uniform (throughout the pathway) steady states and their stability to establish baseline behavior. Section 3 explores stationary distributed (throughout the pathway) profiles, focusing on the influence of activation, inhibition, and bias on steady-state transmission. Section 4 investigates traveling-wave solutions and when they are generated as well as the details of the numerical methods used to estimate propagation speed. Section 5 analyzes the effects of systematic parameter variation and pathway heterogeneity, introducing a velocity-based rescaling to isolate intrinsic propagation. Finally, Section 6 discusses the broader implications of these findings for simplifying biochemical models and outlines future research directions.

\section{Methods}
\label{sec:Methods}
\subsection{Mathematical model}


Biological signaling pathways are fundamentally reaction chains where each node represents the state of a molecular component and each edge defines the regulatory influence of one component on the next. While signaling networks are often characterized by complex cross-talk and feedback loops, this study focuses on individual pathways to isolate the intrinsic dynamics of signal propagation. Building on our previous work \cite{jayathilaka2024two}, we utilize the canonical pathway representation. This framework orients the definition of nodes such that every edge acts as an activator of the subsequent node, regardless of whether the underlying biological relationship is activating or inhibitory. By redefining specific nodes by their ``complement" (e.g., inactive concentration), we treat inhibition as the activation of a node’s inactive form, allowing us to investigate how quickly information in the form of sequential activation traverses the cascade.

In this pathway model, we assign each node a variable $x_i$ for $i = 1, \dots, N$, where $N$ denotes the pathway length. We assume $N$ is sufficiently large to investigate asymptotic propagation rates. Similar to the study of Fisher--KPP waves in infinite media \cite{fisher1937wave}, considering large $N$ allows us to characterize intrinsic wave speeds in the absence of boundary effects, cross-talk, or feedback interference. While these state variables are often mapped to molecular concentrations, we adopt the dimensionless notation $x_i \in [-1, 1]$. This domain emphasizes the functional balance between two opposing molecular forms such as phosphorylated and dephosphorylated species where $+1$ and $-1$ represent the saturation of the signaling component in either respective form. This choice reflects the reality that biochemical information is often encoded by the shift between these binary states. By centering the state at zero, we more effectively characterize how the bias parameter $\phi$ mediates the contribution of each form to the overall signal transduction. For convenience, we retain the nomenclature that positive states constitute activating effects and negative states constitute inhibitory effects; however, in this context, these effects are associated with increases and decreases in $x_i$, respectively. The canonical form ensures that the signal ``direction" is always toward a positive correlation and alignment with the state of upstream nodes.

The evolution of each node is governed by activating inputs from the upstream node $x_{i-1}$. Each edge transition is assumed to operate through a combination of activation (driven by a positive upstream state) and inactivation (driven by a negative upstream state). The relative influence of these two opposing states is controlled by a bias parameter, $\phi_i \in [-1, 1]$. When $\phi = 1$, signal propagation is driven exclusively by the ``on" ($+1$) state of the upstream node; when $\phi = -1$, it is driven by the ``off" ($-1$) state; and $\phi = 0$ represents an unbiased, balanced interaction where both forms contribute equally to the downstream transition. The kinetics are modeled using a reparameterized Michaelis--Menten form, adapted for the $[-1, 1]$ domain. The parameter $\alpha_i$ defines the interaction timescale, while $\beta_i$ represents the saturation characteristics. As detailed in \cite{jayathilaka2024two}, $\beta$ is closely related to the Michaelis constant; $\beta = 1$ yields a linear, uniform response, whereas $\beta \to \infty$ approaches mass-action kinetics.

Quantitatively, the rate of change for a downstream node $x_i$ is given by:

\begin{equation}\label{node_general}\frac{\mathrm{d} x_i}{\mathrm{d} t} = f_{i}(x_{i-1}, x_i) = \frac{(1 + \phi_i)}{4}(1 + x_{i-1})\frac{\alpha_i \beta_i (1 - x_i)}{2\beta_i - (1 + x_i)} - \frac{(1 - \phi_i)}{4}(1 - x_{i-1})\frac{\alpha_i \beta_i (1 + x_i)}{2\beta_i - (1 - x_i)},\end{equation}

where the subscript $i$ indicates parameters specific to a given edge. To initiate signal propagation, we apply an input stimulus at the boundary node $x_0(t)$, the specific form of which is defined in subsequent sections where stationary and traveling-wave profiles are analyzed. All simulations are performed using MATLAB’s ode45 solver. Our analysis begins with the homogeneous case, where all edges share identical parameters ($\phi_i \equiv \phi, \alpha_i \equiv \alpha, \beta_i \equiv \beta$), before we generalize the framework to investigate the impact of edge-to-edge heterogeneity on signal fidelity and speed.

\subsection{Uniform steady states}
\label{sec:steadystate}

Prior to examining traveling waves or parameter-driven changes in propagation speed, it is essential to characterize the fundamental steady-state behavior of the model under uniform parameter conditions. The manner in which nodes are balanced at rest significantly influences the pathway’s response to a stimulus and defines the initial conditions of the system. Furthermore, these steady states describe the equilibrium activity of nodes far from external input signals. Consequently, we begin by considering a uniform regime in which each node stabilizes at an identical activity level. To gain insight into this baseline activity, we assume a pathway of sufficient length such that boundary effects are negligible.

By setting $x_i = x_{i-1} = x$ (with $\phi_i \equiv \phi$, $\alpha_i \equiv \alpha > 0$, and $\beta_i \equiv \beta > 1$), the system reduces to a single autonomous differential equation. This reduction facilitates a clear understanding of the internal pathway pressure prior to the introduction of any heterogeneous effects. Under these uniform conditions, Equation (\ref{node_general}) simplifies to:

\begin{equation}\label{USS_ODE}
\frac{dx}{dt} = \frac{\alpha \beta (1 - x^2) [x + \phi(2\beta - 1)]}{2 [ (2\beta - 1)^2 - x^2 ]}
\end{equation}

Setting $\frac{dx}{dt} = 0$ yields three candidate equilibria:
\begin{equation}
x^*_1 = 1, \qquad x^*_2 = -1, \qquad x^*_3 = -\phi(2\beta - 1) = \xi.
\end{equation}
The interior equilibrium $x^*_3 = \xi$ is biologically meaningful only when it lies within the domain $(-1, 1)$, which requires:
\begin{equation}
-\phi_c < \phi < \phi_c, \quad \text{where} \quad \phi_c = \frac{1}{2\beta - 1}.
\end{equation}

\begin{figure}[H]
\centering
\includegraphics[width=0.8\textwidth]{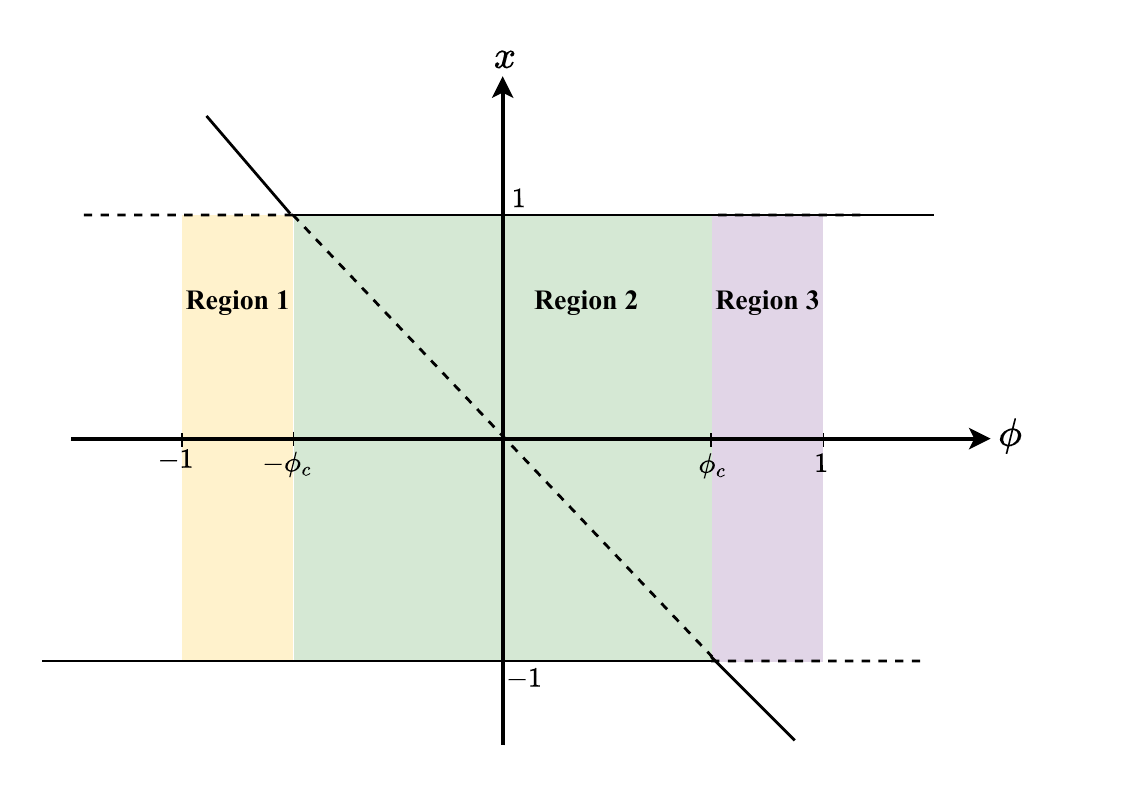}
\caption{Bifurcation diagram showing uniform steady states as a function of pathway bias $\phi$. Stable equilibria are indicated by solid lines and unstable equilibria by dashed lines. In the extreme bias regimes, Region 1 ($\phi < -\phi_c$) and Region 3 ($\phi > \phi_c$), the pathway is monostable at $x=-1$ and $x=1$, respectively. Region 2 ($-\phi_c < \phi < \phi_c$) displays bistability, where the unstable equilibrium $x = \xi$ (Equation \ref{USS_ODE}) serves as the threshold between the two stable basins.}\label{fig:pplane}
\end{figure}

Figure \ref{fig:pplane} illustrates the bifurcation of these uniform steady states as a function of the bias parameter $\phi$. The system exhibits three distinct dynamical regimes. \textbf{Region 1 ($\phi < -\phi_c$)} and \textbf{Region 3 ($\phi > \phi_c$)} represent extreme bias regimes where the system is monostable; the states $x = -1$ and $x = 1$ are globally attracting, respectively. In these regions, the pathway is effectively ``locked,'' and a signal must overcome significant internal bias to propagate. 

\textbf{Region 2 ($-\phi_c < \phi < \phi_c$)} is defined by the coexistence of all three equilibria. In this regime, the system exhibits bistability between the fully active state ($x=1$) and the fully inactive state ($x=-1$). The interior equilibrium $x=\xi$ acts as an unstable steady state (a separatrix), defining the threshold that an upstream signal must exceed to flip a downstream node from one basin of attraction to the other. As $\phi$ crosses the critical values $\pm \phi_c$, the system undergoes transcritical bifurcations where the unstable equilibrium $\xi$ collides with one of the boundary steady states, exchanging stability and fundamentally altering the pathway's capacity for signal transmission.

\subsection{Stationary distributions}
\label{sec:statdist}

The uniform steady states described in Section \ref{sec:steadystate} characterize the limiting behavior of the system as both $t \to \infty$ and $i \to \infty$; they represent the activity level of nodes sufficiently distant from the pathway boundary to be unaffected by the input stimulus. However, in the presence of a constant external stimulus at the first node ($i=1$), the pathway develops a stable \textit{stationary spatial distribution}, a node-dependent profile that describes the steady-state activation levels across the cascade at long times. We initially examine these distributions under uniform pathway conditions ($\alpha_i \equiv \alpha, \beta_i \equiv \beta, \phi_i \equiv \phi$). These stationary profiles are of fundamental importance to our subsequent analysis of signal propagation, as they define the state of the pathway in the wake of an advancing traveling wave. While our primary focus remains on canonical pathways where nodes act as sequential activators \cite{jayathilaka2024two}, it is essential to characterize the pathway's response to both activating and inhibitory stimuli to fully understand the dynamical range of the oriented form.

The input stimulus at the first node ($i=1$) is modeled based on the framework established in \cite{jayathilaka2024two}. The rate of change for $x_1$ is determined by Equation (\ref{node_general}), where the upstream state $x_0$ represents a constant external signal. We define activation by $0 < x_0 \leq 1$ and inhibition by $-1 \leq x_0 < 0$, where the magnitude of $x_0$ determines the strength of the input relative to the intrinsic catalytic rate $\alpha$.

\begin{equation}\label{node_initial}
\frac{\mathrm{d} x_1}{\mathrm{d} t} = \frac{(1 + \phi_1)}{4}(1 + x_{0})\frac{\alpha_1 \beta_1 (1 - x_1)}{2\beta_1 - (1 + x_1)} - \frac{(1 - \phi_1)}{4}(1 - x_{0})\frac{\alpha_1 \beta_1 (1 + x_1)}{2\beta_1 - (1 + x_1)}.
\end{equation}

When investigating the response of the pathway to a new stimulus, we assume the system initially resides at a uniform ``resting" level of activation. We consider two primary initial condition cases: (1) the fully inactive state, $x_i(0) = -1$, and (2) the fully active state, $x_i(0) = 1$. The choice of initial condition is informed by the bias parameter $\phi$ and the resulting stability regions identified in Fig. \ref{fig:pplane}. In Regions 1 and 3, where $x=-1$ and $x=1$ are globally attracting, cases (1) and (2) respectively represent the only stable ``no-input" configurations. In the bistable Region 2, either initial condition may be used; here, the pathway's response depends on which stable equilibrium the system originates in, allowing us to study how activation waves propagate through an inactive pathway, or how inhibitory signals might trigger a receding wave in an active one.

\begin{figure}[h]
\centering
\includegraphics[width=1\textwidth]{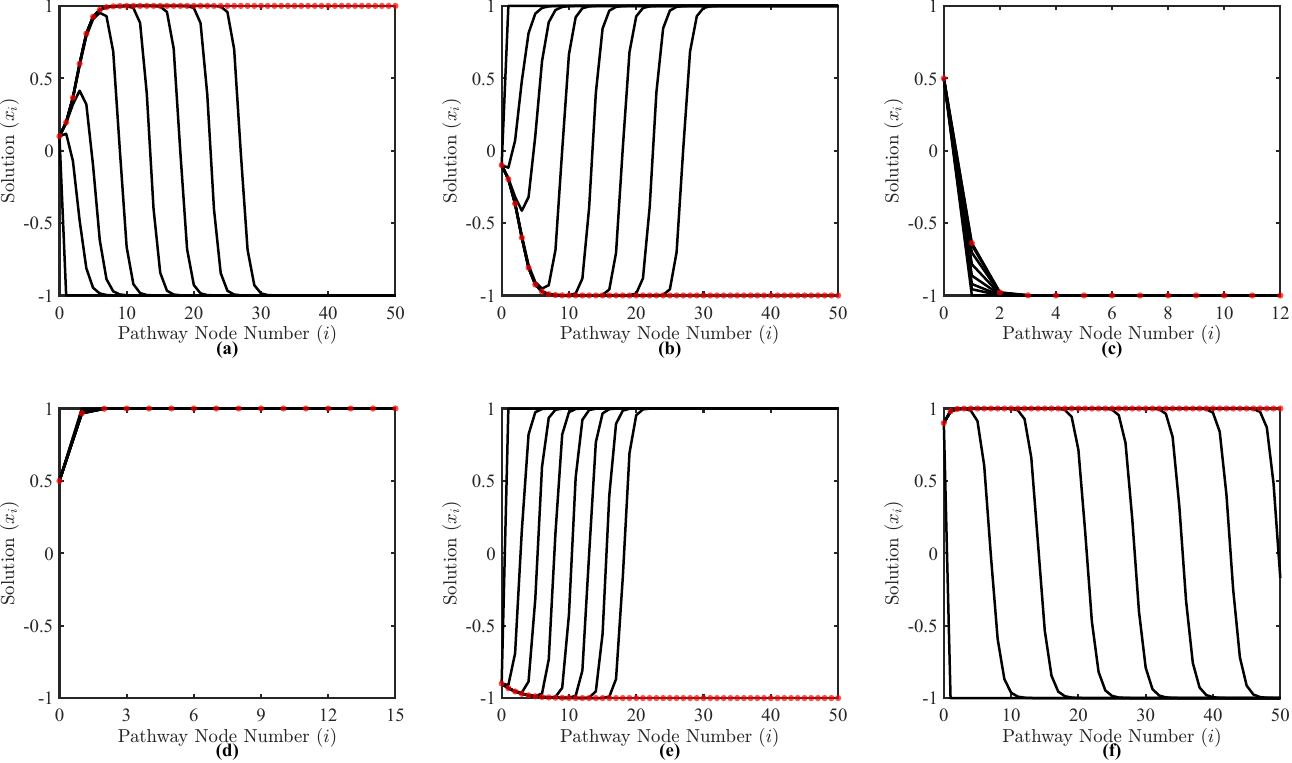}
\caption{Transient and stationary dynamics of node activity $x_i$ across a uniform pathway ($N=200, \alpha=1, \beta=1.5, \phi_c=0.5$). Black curves represent successive temporal snapshots of $x_i(t)$, illustrating the evolution of the spatial profile toward a long-time attractor. Red markers denote the analytical stationary distribution computed iteratively from Eq.~\eqref{quadraticSol2}. 
(a--b) Unbiased case ($\phi=0$, Region 2): Symmetric activation and inhibition stimuli produce smooth, gradually evolving fronts that converge algebraically toward the active and inactive states, respectively.
(c--d) Extreme bias regimes: In Region 1 ($\phi=-0.75$), negative bias rapidly suppresses activation; in Region 3 ($\phi=+0.75$), positive bias drives near-instantaneous alignment with the active state. Both regimes lack sustained wave propagation, relaxing locally to the globally stable equilibrium.
(e--f) Near-critical bias ($\phi=0.3$, Region 2): Symmetry-breaking results in a slower inhibitory front in (e) and a faster activatory front in (f). 
In all cases, the close agreement between the late-time numerical snapshots and the red markers validates the iterative map $g(x_{i-1})$ as the definitive spatial attractor for the signaling cascade.}
\label{fig:stationary_sol}
\end{figure}

Figure \ref{fig:stationary_sol} illustrates the transient dynamics of node activity $x_i(t)$ in a uniform pathway of $N=200$ nodes ($\alpha=1, \beta=1.5, \phi_c=0.5$). Across all dynamical regions, the temporal snapshots (black curves) converge toward the analytical stationary distribution (red markers) computed iteratively using Eq. \eqref{quadraticSol2} derived later in this section, confirming it as the long-term spatial attractor.

In the unbiased case ($\phi=0$, Region 2), panels (a) and (b) show that activation and inhibition stimuli both produce smooth, symmetric traveling wavefronts that propagate at constant speeds. In the wake of these fronts, the pathway converges slowly to the stationary distribution, reflecting the algebraic convergence characteristic of the unbiased regime. Conversely, in the extreme bias regimes (Region 1, panel c; Region 3, panel d), the system does not support sustained wave propagation. Instead, the strong internal bias causes the signal to decay or relax locally and rapidly to the globally stable equilibrium.

Panels (e) and (f) demonstrate the effect of a near-critical positive bias ($\phi=0.3$) in Region 2. Here, the symmetry between activation and inhibition is broken: the positive bias reinforces the activatory input, resulting in a faster wavefront in (f). The same bias slows down the transmission speed in the case of inhibitory input in (e) where the magnitude of inhibition at the input necessarily has to overcome the significantly negative separatrix for this level of positive bias. At long times, the agreement between the numerical profiles in red from  of Eq. \eqref{quadraticSol2} and time dependent ones in black  validates our analytical framework for predicting the pathway's spatial steady state across diverse kinetic regimes which follows in this section.

We can find approximate stationary pathway distributions in the limit as $N\rightarrow\infty$. Consider the dynamic model for node $i>1$ under uniform pathway conditions according to Equation (\ref{node_general}):

\begin{equation}
    \frac{d x_i}{d t} = \phi^{+} \frac{(1+x_{i-1})(1-{x_i})}{2 \beta-(1+x_i)}-\phi^{-} \frac{(1-x_{i-1})(1+x_i)}{2 \beta-(1-x_i)},
\label{eq:inhibition_repeat}
\end{equation}

where $\phi^{\pm} = \alpha \beta (1\pm\phi) / 4$ for notational simplicity. Denote $B=2 \beta-1 > 1$ and $\Phi = \phi^{+}/\phi^{-} > 0$. Solving $\frac{d x_i}{d t} = 0$ reduces Equation (\ref{eq:inhibition_repeat}) to:

\begin{equation} \label{reducedStationaryBalanceGeneral}
   \Phi (B+x_{i})(1+x_{i-1})(1-x_{i})=(B-x_{i}) (1-x_{i-1}) (1+x_{i}).
\end{equation}

We can define the iterative equation for the stationary distribution by solving the quadratic in $x_i$ in terms of $x_{i-1}$. Specifically, it is advantageous to define $X^\pm_{i-1} = \Phi(1+x_{i-1}) \pm (1-x_{i-1})$ and let $\chi_{i-1} = X^+_{i-1}/X^-_{i-1}$, which simplifies Equation (\ref{reducedStationaryBalanceGeneral}) into the more manageable form:

\begin{equation}\label{quadraticsimplified} 
x_i^2 +(B-1)\chi_{i-1} x_i - B = 0.
\end{equation}

Since $B>1$ and $|\chi_{i-1}| > 1$ (a consequence of $\Phi >0$ and $|x_{i-1}|<1$), Equation (\ref{quadraticsimplified}) always has one and only one root satisfying $|x_i| < 1$. To find this root, Vieta's formulas and the Intermediate Value Theorem indicate that the larger root is chosen if $x_{i-1} > -\phi$ ($\chi_{i-1}>0$), otherwise the smaller root is chosen. We embed this choice analytically by taking the positive branch but keeping $\chi_{i-1}$ factored outside the square root of the discriminant:

\begin{align}
   x_i &= \frac{\chi_{i-1}}{2}\left(-(B-1) + \sqrt{(B-1)^2 + \frac{4B}{\chi_{i-1}^2}}\right),\label{quadraticSol}\\ 
   &= \frac{1}{2(\phi+x_{i-1})}\left(-(B-1)(1+\phi x_{i-1}) + \sqrt{(B-1)^2 (1+\phi x_{i-1})^2 + 4B(\phi+x_{i-1})^2}\right) \label{quadraticSol2} \\ 
   &= g(x_{i-1};\phi,B).
\end{align}

Importantly, the mapping $g(x;\phi,B)$, as defined by (\ref{quadraticSol2}), satisfies $g(x;\phi,B) < x$ for all $-1 < x < \xi$ and $g(x;\phi,B) > x$ for all $\xi < x < 1$. This guarantees that the stationary distribution $x_i$ is monotonically decreasing toward $-1$ if $x_0 < \xi = -\phi B$, and monotonically increasing toward $1$ if $x_0 > \xi$. This asymptotic behavior is entirely consistent with the uniform steady-state bifurcation diagram shown in Fig. \ref{fig:pplane}.

To quantify this spatial relaxation toward the final stationary distribution, we first examine the asymptotic behaviour far downstream from the input. As the signal approaches the stable uniform steady state $x^*=\pm1$, the deviation from this equilibrium, $\epsilon_i = |x^* - x_i|$, decays geometrically. The exact rate of this geometric decay, $\lambda$, is dictated by the local linear stability of the fixed point. This is found by evaluating the first derivative of the exact spatial mapping function $g(x; \phi, B)$, defined in (\ref{quadraticSol2}), precisely at the steady state $x = x^*$: $ \lambda(B, \phi) = \left| g'(x^*) \right| $. Using (\ref{quadraticSol2}) we find
\begin{equation} \label{lambdaEq}
\lambda(B, \phi) = \left(\frac{1-\phi}{1+\phi}\right)^{x^*} \left( \frac{B-1}{B+1} \right).
\end{equation}

While this geometric decay exactly characterizes the asymptotic tail far from the input node at $i=1$, we must also examine the initial transient phase, which determines how deeply an input signal effectively penetrates the pathway before asymptotic decay begins. Expanding Equation (\ref{quadraticSol2}) asymptotically in the large coupling limit ($B \to \infty$) yields:
$$ x_i \approx x_{i-1} + \frac{(\phi+x_{i-1})(1-x_{i-1}^2)}{B(1+\phi x_{i-1})} $$
This naturally explains why sigmoid-like distributions and slower convergence rates emerge at large $B$: the discrete update step $\Delta x$ is inversely proportional to $B$. This vanishing step size allows for smooth, extended spatial profiles that linger around the input value $x_0$ deep into the pathway before eventually accelerating into the geometric tail.

To gain insight into this transient stagnation, we reduce our analysis to the perfectly symmetric, unbiased limit ($\phi=0$). In this regime, the separatrix sits exactly at the origin ($\xi=0$). The system converges to the active steady state $x^*=1$ for any $x_0>0$, and behaves symmetrically, converging to $x^*=-1$ for any $x_0<0$. In either case, the convergence rate simplifies to $ \lambda = (B-1)(B+1)^{-1} $. To model the initial lingering phase, we introduce a spatial delay parameter, $\Delta i$. This delay acts as a horizontal shift on the geometric decay, allowing us to express the stationary distribution far from the input as:
$$ \epsilon_i \approx (1-|x_0|)\lambda^{i-\Delta i}. $$

To derive an approximation for $\Delta i$, we appeal to a continuous spatial limit of the unbiased system, governed by the differential equation $\frac{dx}{di} \approx g(x; 0, B) - x = \frac{x(1-x^2)}{B}$ (which notably matches our large-$B$ expansion when $\phi=0$). Separating variables and integrating from the initial state $x_0$ to a downstream state $x_i$ yields the traversal distance in the continuous approximation:
$$ i = B\ln\left(\frac{x_i}{x_0}\sqrt{\frac{1-x_0^2}{1-x_i^2}}\right). $$

As the signal asymptotically approaches the active steady state ($|x_i| \to 1$), we substitute $|x_i| = 1-\epsilon_i$. Recognizing that in this limit $1-x_i^2 \approx 2\epsilon_i$, the integral simplifies to:
$$ i \approx \frac{B}{2}\ln\left(\frac{1-x_0^2}{2x_0^2\epsilon_i}\right). $$

Inverting this expression for $\epsilon_i$ reveals the asymptotic decay strictly within the continuous framework:
$$ \epsilon_i \approx \left(\frac{1-x_0^2}{2x_0^2}\right)e^{-2i/B}. $$

To bridge this continuous approximation back to the discrete reality of the lattice, we map the continuous spatial decay term $e^{-2/B}$ to the exact discrete multiplier $\lambda$. Equating this resulting continuous expression to our spatially delayed geometric ansatz, $(1-|x_0|)\lambda^{i-\Delta i}$, allows the $(1-|x_0|)$ terms to neatly cancel. Solving for the shift yields a closed-form analytical approximation for the spatial penetration depth:
\begin{equation} \Delta i \approx \frac{1}{\ln\lambda}\ln\left(\frac{2x_0^2}{1+|x_0|}\right). \label{Deltai}\end{equation}

This analytical form explicitly demonstrates the physics of the separatrix: as the initial state approaches the unstable fixed point at the origin ($|x_0| \to 0$), the argument of the logarithm diverges proportionally to $1/x_0^2$. This correctly predicts a rapidly growing spatial delay as signals struggle to escape the unstable equilibrium.

\begin{figure}[htbp]
\centering
 \includegraphics[width=\linewidth]{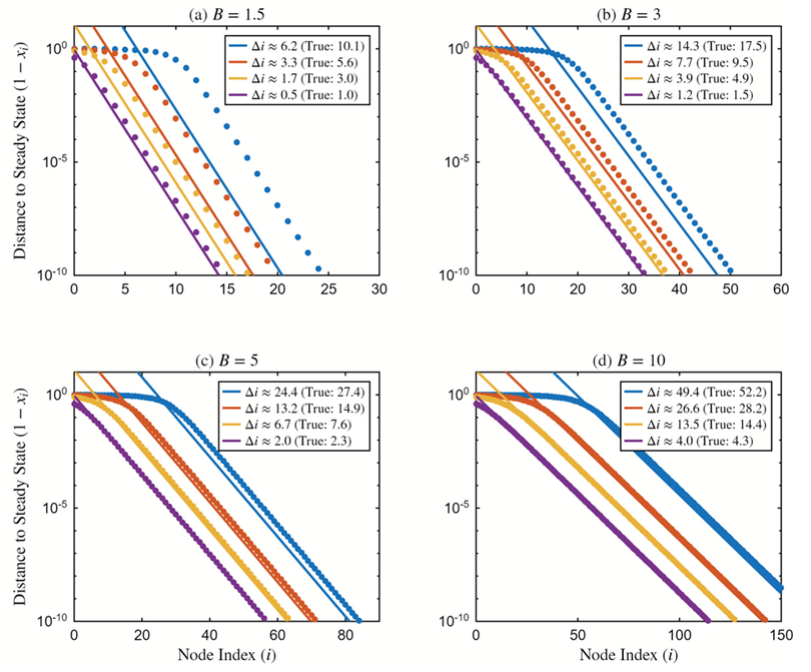}
\caption{Stationary distributions in the unbiased regime ($\phi=0$) for varying coupling parameters: (a) $B=1.5$, (b) $B=3.0$, (c) $B=5.0$, and (d) $B=10.0$. Markers denote the exact numerical evaluation using Equation (\ref{quadraticSol2}) and each subplot contains initial states ranging in distance from the separatrix given by $x=\xi=0$: $x_0 = 0.005$ (blue), $x_0 = 0.05$ (orange), $x_0 = 0.2$ (yellow), and $x_0 = 0.6$ (purple). Solid lines represent the corresponding analytical asymptotic predictions, $\epsilon_i \approx (1-|x_0|)\lambda^{i-\Delta i}$, utilizing the exact discrete convergence rate (\ref{lambdaEq}) and approximate penetration depth $\Delta i$ using Equation (\ref{Deltai}).  The value of the approximate and true penetration depths are printed in the legend. The true penetration depth is extracted via tail-fitting the exact stationary distribution whilst assuming the asymptotic rate $\lambda$. The approximate penetration depth $\Delta i$ scales appropriately but is more inaccurate as $B$ approaches 1 and $x_0$ approaches the separatrix $\xi$. The approximate penetration depth consistently provides a lower bound on the true penetration depth}
\label{fig:spatial_delay}
\end{figure}

To validate this analytical scaling and compare our approximations against numerically exact stationary distributions, we iteratively evaluate the discrete mapping (Equation \ref{quadraticSol2}) for $\phi=0$ across a range of initial conditions $x_0>0$ and saturation parameters $B>1$. Fig. \ref{fig:spatial_delay} illustrates the exact spatial decay (markers) alongside the analytical predictions incorporating $\Delta i$ (solid lines).

The visual evidence highlights two critical system dynamics. First, there is perfect agreement between the numerical mapping and the analytical geometric rate of convergence (Equation \ref{lambdaEq}). Second, the analytical penetration shift $\Delta i$ (Equation \ref{Deltai}) successfully captures the overarching scaling with respect to both $B$ and $x_0$. While the continuous approximation slightly underestimates the absolute shift as $x_0$ approaches the separatrix $\xi$ (an expected artifact due to the step-wise stagnation of the discrete map), the relative error remains tightly bounded.

The stationary distributions characterized above describe the long-time spatial profile that develops in the wake of a traveling front---the state toward which the pathway ultimately converges following the introduction of a stimulus. However, they provide no information regarding the \emph{rate} at which this profile is established, nor how quickly the signaling front itself propagates through the cascade (if it cascades at all---see Fig. \ref{fig:stationary_sol}). Because the time taken for a signal to traverse the pathway determines the speed and reliability of biological information transfer, we turn now to the question of wave propagation velocity, examining its dependence on the kinetic parameters $\alpha$, $\beta$, and $\phi$.

\subsection{Signal propagation speed in uniform pathways}
\label{sec:wavespeed}

While the preceding sections have characterized the long-time stationary spatial profile of a pathway under a sustained input, a critical functional property of biological signalling is the dynamic time delay associated with finite signal propagation following a localized stimulus. In this section, we investigate the velocity of this propagating wave as a function of the underlying chemical kinetic parameters, restricting our attention to the uniform pathway ($\phi_i \equiv \phi$, $\alpha_i \equiv \alpha$, $\beta_i \equiv \beta$).

We restrict our dynamic analysis to Region~2 ($|\phi| < \phi_c$), where the two stable uniform steady states $x^* = -1$ and $x^* = +1$ coexist. As demonstrated previously in Figure~\ref{fig:stationary_sol}(c) and (d), Regions~1 and~3 do not support sustained travelling waves; because only a single stable uniform state exists in these regimes, the pathway relaxes rapidly and locally to the globally attracting equilibrium without forming a propagating front. In Region~2, by contrast, a well-defined travelling wavefront can exist, connecting the two stable states and advancing downstream at a measurable speed. 

The scenario of interest is one in which the pathway is initialized entirely in one stable uniform steady state, and a sustained boundary input $x_0$ is applied that lies in the basin of attraction of the complementary steady state. Furthermore, as we are focused on signal wave transmission speed in this section we will further restrict the input $x_0$ to $\pm 1$ so that the wavefront is not distorted by a non-uniform stationary distribution left in the wake of the wave. Under these conditions, a travelling wavefront propagates downstream, sequentially switching nodes from the initial state $x_i(0) = \pm 1$ toward the final active or inactive state, with the stationary distribution described in Section~\ref{sec:statdist} developing in the wake of the front. Our objective here is to analytically and numerically quantify the propagation speed of this front as a function of the underlying kinetic parameters.


In an idealized infinite pathway, a travelling wave solution would naturally assume the form $x_i(t) = f(i - ct)$, where $c$ denotes the constant signal propagation velocity (in nodes per unit time). However, in our finite simulated cascade of length $N$, we do not expect a perfect, shape-invariant travelling wave, particularly near the pathway boundaries where boundary conditions and transient stationary distribution exist. 

To quantify the wave speed, we estimate an instantaneous velocity $c_j$ at discrete time points $t_j$ (separated by a constant step $\delta t = t_{j+1} - t_j$). We achieve this by finding the spatial shift $c \delta t$ that minimises the least-squares difference between the pathway profile at $t_{j+1}$ and the shifted profile from $t_j$. Specifically, we define the instantaneous velocity as:
$$ c_j = \arg \min_{c > 0} F_j(c), $$
where the objective function $F_j(c)$ is given by:
\begin{equation}
    F_j(c) = \sum_{i=1}^{N} \left( x_i(t_{j+1}) - \tilde{x}(i - c\,\delta t, t_j) \right)^2.
    \label{eq:cj}
\end{equation}
Here, $\tilde{x}(\cdot, t_j)$ denotes a linear spatial interpolation of the discrete pathway profile at time $t_j$, necessitated by the fact that the upstream shifted coordinate $i - c\,\delta t$ is generally non-integer. 

By solving this optimisation numerically at each time step, we extract a time series of instantaneous velocity estimates, $\{c_j\}$, that captures the dynamic evolution of the propagating front. This approach leverages the approximate shape-invariance of the travelling wave while remaining robust to the deformations that occur near the cascade extremities.

\begin{figure}[h]
    \centering
    \subfigure[]{%
        \includegraphics[width=0.52\textwidth]{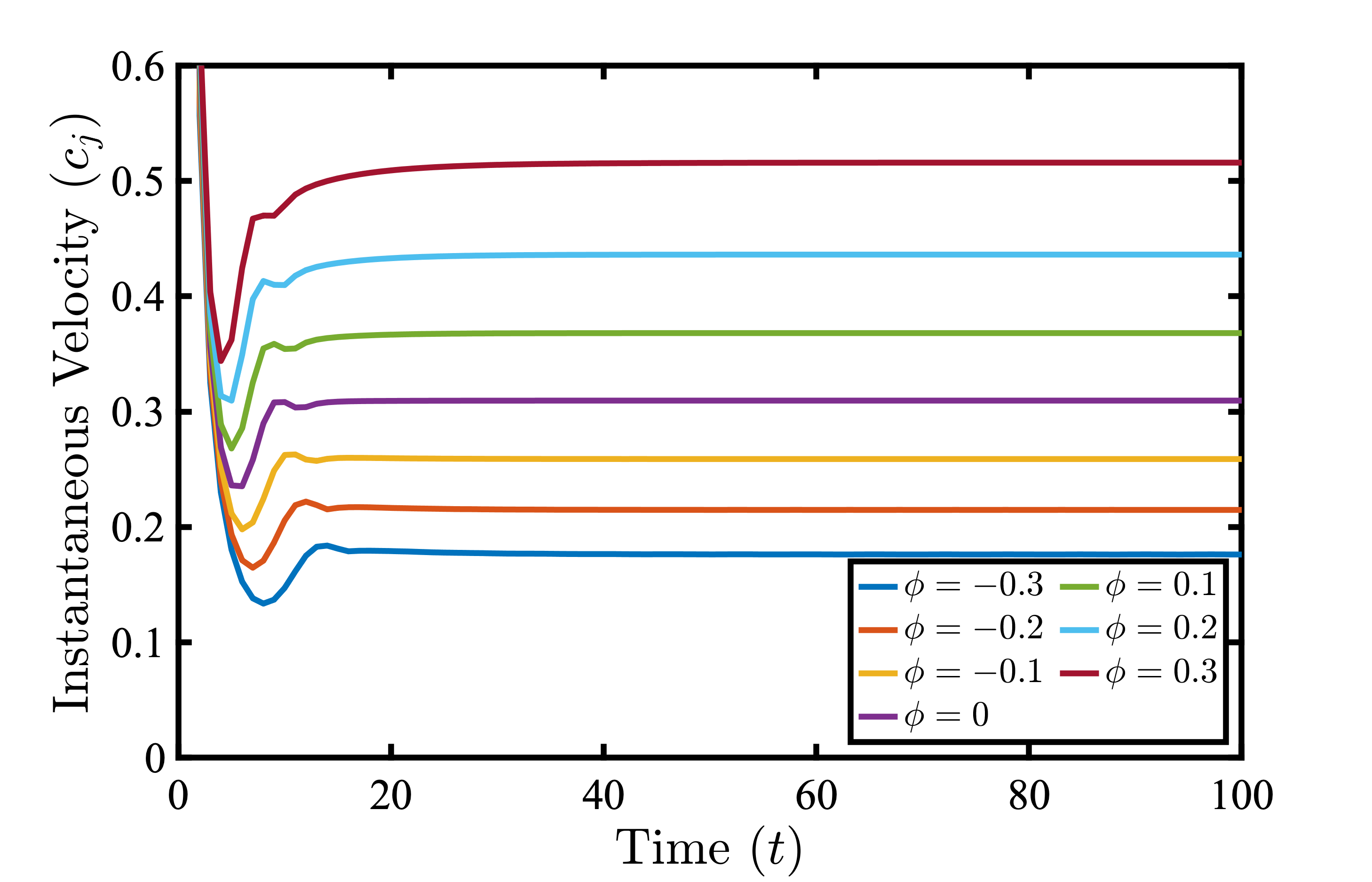}
        \label{fig:vel_transient}
    }%
    \hfill
    \subfigure[]{%
        \includegraphics[width=0.45\textwidth]{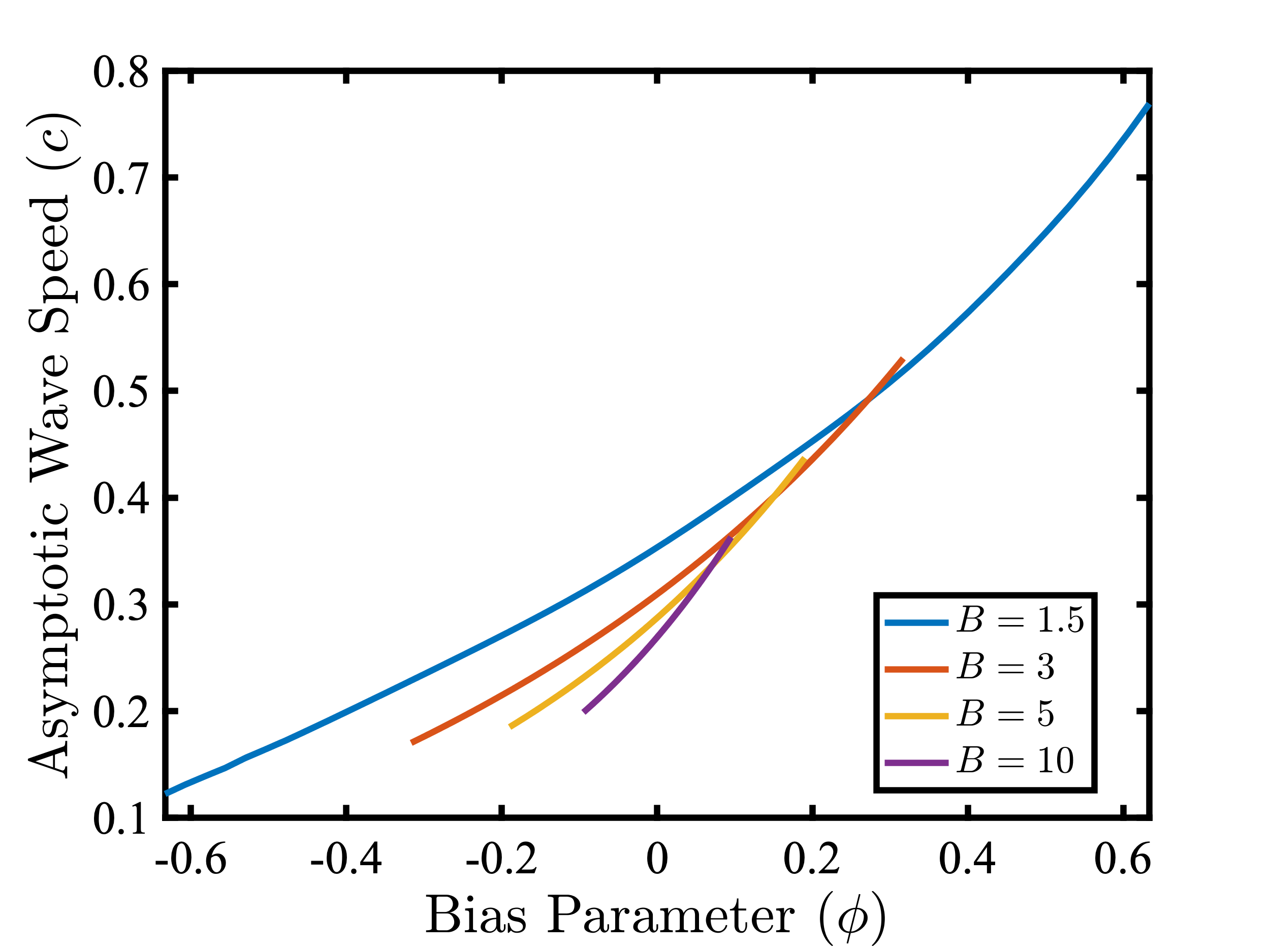}
        \label{fig:vel_asymptotic}
    }
\caption{\textbf{Transient and asymptotic wave speeds of signal propagation in the bistable regime.} (a) Instantaneous wave speed over the initial $100$ time units for a fixed nonlinearity parameter $B = 3$. (b) Asymptotic wave speed as a function of the bias parameter $\phi$ for multiple nonlinearity values $B \in \{1.5, 3, 5, 10\}$. To eliminate boundary distortions, the asymptotic wave speeds are evaluated at exactly $50\%$ of the wave's total propagation time. The wave speed is only considered between the critical bifurcation boundaries of Region 2 where bistability of the uniform steady state is a key feature $\pm \phi_c = \pm 1/B$. For all simulations, the cascade length is $N = 200$ nodes, the activation rate is $\alpha = 1$, the initial uniform state is $x_i(0) = -1$, and the driving boundary input is $x_0 = 1$.}
\label{fig:velocity}
\end{figure}

We apply this velocity tracking method to simulated cascades within the bistable regime (Regime 2) to characterize both the transient and steady-state propagation dynamics. For all simulations in this section, we fix the cascade length to $N = 200$ nodes and the activation rate to $\alpha = 1$, with the network uniformly initialized at $x_i(0) = -1$ and subjected to a sustained boundary input of $x_0 = 1$. Figure~\ref{fig:velocity}(a) illustrates the instantaneous wave speed over the initial $100$ time units for a fixed nonlinearity parameter $B = 3$ across various bias parameters $\phi$. As the wavefront develops and travels downstream, the propagation speed quickly transitions from an initial boundary-dominated transient phase into a constant asymptotic wave speed. To accurately reflect the physical limits of the finite cascade, the instantaneous velocity tracking for each $\phi$ is actively halted the moment the wavefront reaches the terminal node (defined as a deviation of $10^{-4}$ from its initial state).

Building on this, Figure~\ref{fig:velocity}(b) quantifies the final asymptotic wave speed as a function of the bias parameter $\phi$ for multiple nonlinearity values $B \in \{1.5, 3, 5, 10\}$. For each curve, $\phi$ is varied strictly within its corresponding bistable window, bounded by $\pm \phi_c$, where $\phi_c = 1/B$. To ensure the speed is measured after the transient phase has fully decayed but well before boundary effects from the terminal node distort the wave, the stationary speed is extracted at exactly $50\%$ of the total propagation time---the temporal midpoint of the wave's journey through the cascade.

The asymptotic wave speed profiles presented in Figure~\ref{fig:velocity}(a) reveal several key physical insights into the signal transmission dynamics. For the activation wave under consideration driven by a boundary input of $x_0 = 1$ into a cascade initialized at $x_i(0) = -1$---the propagation speed monotonically increases as the bias parameter $\phi$ is varied from $-\phi_c$ to $\phi_c$. Physically, as the wavefront reaches a given node, the active upstream neighbor must pull the local state from $-1$ to $1$ against the system's natural restorative bias toward the inactive steady state. Crucially, linear stability analysis of the leading edge reveals that the resting state at $x = -1$ is strongly stable against small perturbations. Consequently, the cascade exclusively supports pushed waves, where propagation is driven not by a marginally stable leading edge, but by the nonlinear bulk of the wave brute-forcing its way across the intermediate unstable steady state, $\xi$. For smaller values of $\phi$, this intermediate state $\xi$ is shifted toward the active state, making the restorative interval $x_i \in [-1, \xi]$ wider. This presents a larger barrier for the wavefront bulk to overcome, thereby attenuating the propagation speed. Conversely, as $\phi$ increases, this interval narrows, facilitating much faster wave transmission.

Due to the inherent symmetry of the bistable kinetic formulation, it is unnecessary to independently analyze the complementary inhibitory wave (where $x_0 = -1$ propagates into a network initialized at $x_i(0) = 1$). The speed of an inhibitory wave at a given bias $\phi$ is strictly equivalent to the speed of an activation wave at $-\phi$. Therefore, the parameter sweep in Figure~\ref{fig:velocity}(b) provides a complete characterization of both activation and inhibition signal velocities.

The nonlinearity parameter $B$ dictates both the global scaling of the propagation speed and the spatial morphology of the wavefront. Increasing $B$ generally decelerates the wave across the parameter space. This behavior reflects the heightened saturation of the underlying chemical kinetics, which effectively dampens the overall reaction rates of the network. Analytically, higher values of $B$ steepen the transition between the stable roots, heavily penalizing intermediate states and resulting in a severely compressed wavefront spatial width. Because the speed of a pushed wave depends integrally on the nonlinear transition region, this spatial narrowing restricts the domain over which the driving gradient operates, slowing down the signal. However, this kinetic sluggishness can be partially overridden if the system is sufficiently biased. As $\phi$ approaches the upper critical boundary $\phi_c$, the intrinsic tendency toward the positive steady state becomes overwhelming; this structural bias compensates for the kinetic saturation, allowing the signal to propagate rapidly even for large values of $B$.

Finally, while $B$ and $\phi$ strictly govern the nonlinear shape, spatial width, and relative driving forces of the wavefront, the parameter $\alpha$ acts as a purely temporal scaling factor. Because $\alpha$ uniformly scales the kinetic rates governing both the driving and restorative forces, it directly dictates the fundamental time scale of the system. Consequently, the asymptotic wave speed scales linearly with $\alpha$, leaving the spatial wavefront width and the critical threshold boundaries entirely invariant.

\section{Signal propagation speed in heterogeneous pathways}
To systematically analyze signal transmission in heterogeneous networks, we introduce an \textit{a priori} spatial rescaling transformation designed to normalize the pathway. This approach relies on a significant heuristic assumption: despite the complex, multi-node interactions that comprise the nonlinear wavefront, the wavefronts themselves are spatially very narrow typically spanning only a single-digit number of nodes. Because the transition region is so thin, we assert that the wave speed associated with a heterogeneous pathway \textit{should} be intrinsically approximated at the level of individual edges. Specifically, we define a theoretical speed associated with a given edge $i$ in a heterogeneous pathway, $c(\alpha_i, B_i, \phi_i)$, which is the velocity a wave would achieve in a strictly homogeneous pathway where all edges have the same parameters as the specific edge $i$ in a heterogeneous pathway.

We therefore assert that if we define an internode pathway ``distance'' $\Delta s_i = s_{i+1} - s_i$ to be inversely proportional to the theoretical intrinsic connecting edge speed, information will travel down a heterogeneous pathway in a more uniform way. This recursively generates a rescaled position coordinate $s_i$ for each node:
\begin{equation}
\label{eq:rescaling}
    s_{i+1} = s_i + \Delta s_i = s_i + \frac{\bar{c}}{c(\alpha_i, B_i, \phi_i)},
\end{equation}
where $s_0 = 0$, and $\bar{c}$ is a fixed increment chosen to recast the original discrete spatial domain $i \in [0, N]$ onto a discrete set on the unit interval $s \in [0, 1]$. Specifically, setting $\bar{c} = \left( \sum_{k=0}^{N-1} c(\alpha_k, B_k, \phi_k)^{-1} \right)^{-1}$ ensures the total rescaled pathway length is strictly well-defined and equal to unity.

By enforcing this inverse proportionality, the transformation stretches the spatial coordinate where the local kinetics are slow (small theoretical $c$) and compresses it where they are fast (large theoretical $c$). Meanwhile, let $c_j$ denote the instantaneous wave velocity estimated at time $t_j$. Our central hypothesis is that plotting the heterogeneous pathway state $x_i(t)$ against the rescaled coordinate $s_i$, rather than the discrete node index $i$, will collapse the non-uniform signal onto a signal that moves with approximately a fixed velocity whilst not significantly distorting the waveform.

To quantitatively measure these properties, we must formally define the instantaneous wave velocity and the metric for shape variation.

We compute the instantaneous speed $c_j$ in the original pathway using Equation (\ref{eq:cj}). In the normalised pathway the normalised speed $\tilde{c}_j$ at the same timestep $j$ is computed using the definition  we require a modification to this definition
$$ \tilde{c}_j = \arg \min_{c > 0} \tilde{F}_j(c), $$
where the objective function $\tilde{F}_j(c)$ is given by
\begin{equation}
    \tilde{F}_j(c) = \sum_{i=1}^{N} \left( x_i(t_{j+1}) - \tilde{x}^{(s)}(s_i - c\,\delta t, t_j) \right)^2,
    \label{eq:cjtilde}
\end{equation}
where $\tilde{x}^{(s)}(s, t_j)$ denotes a linear spatial interpolation of the discrete pathway profile at time $t_j$ and on the rescaled positions defined by $\{s_i\}_{i=1}^N$. To make a fair comparison between instantaneous speed in the original pathway and the normalized pathway, we compare $c_j/N$ to $\tilde{c}_j$ since the pathway spatial domain has been mapped from an interval of length $N$ to one of unit length.

To evaluate the structural integrity of the wavefront over the entire domain, we introduce a  Global Shape Residual  relative to a wavefront at a particular moment in time. As we expect a more uniform travelling wave in the normalized pathway, We define the moment $t_J$ to be the time in which  $\tilde{x}^{(s)}(0.5, t_J) = 0$ (that is, when the center of the wave passes over the center of the normalized pathway). We then define the Global Shape Residual of the wave at time $t_j$ relative to the wavefront at time $t_J$ (in the figures in the proceeding subsections, this reference wave is indicated by a black dashed line) as the minimum square sum of state differences after translation of the original pathway at time $t_j$ to the wave at time $t_J$ and weighted by the uniform relative distances for each node $1/N$ and call this $R_j$. For the normalized pathway, the Global Shape Residual $\tilde{R}_j$ is computed in the same way but translating in the normalized domain and weighted by $\Delta s_i$ term. 
\begin{align}
\label{eq:shape_res}
    R_j &= \min_{\delta} \sum_{1}^{N} \left( x_i(t_{j}) - \tilde{x}(i - \delta, t_j) \right)^2/N, \\
    \tilde{R}_j &= \min_{\delta} \sum_{1}^{N} \left( x_i(t_{j}) - \tilde{x}^{(s)}(s_i - \delta, t_j) \right)^2 \Delta s_i.
\end{align}
If the signal propagates as a perfect traveling wave with absolute shape invariance, the shifted profile matches exactly, yielding $R_j = 0$.

\subsection{Pathways with Monotonic Parameter Variation}

\begin{figure}[h]
    \centering
    \includegraphics[width=1\textwidth]{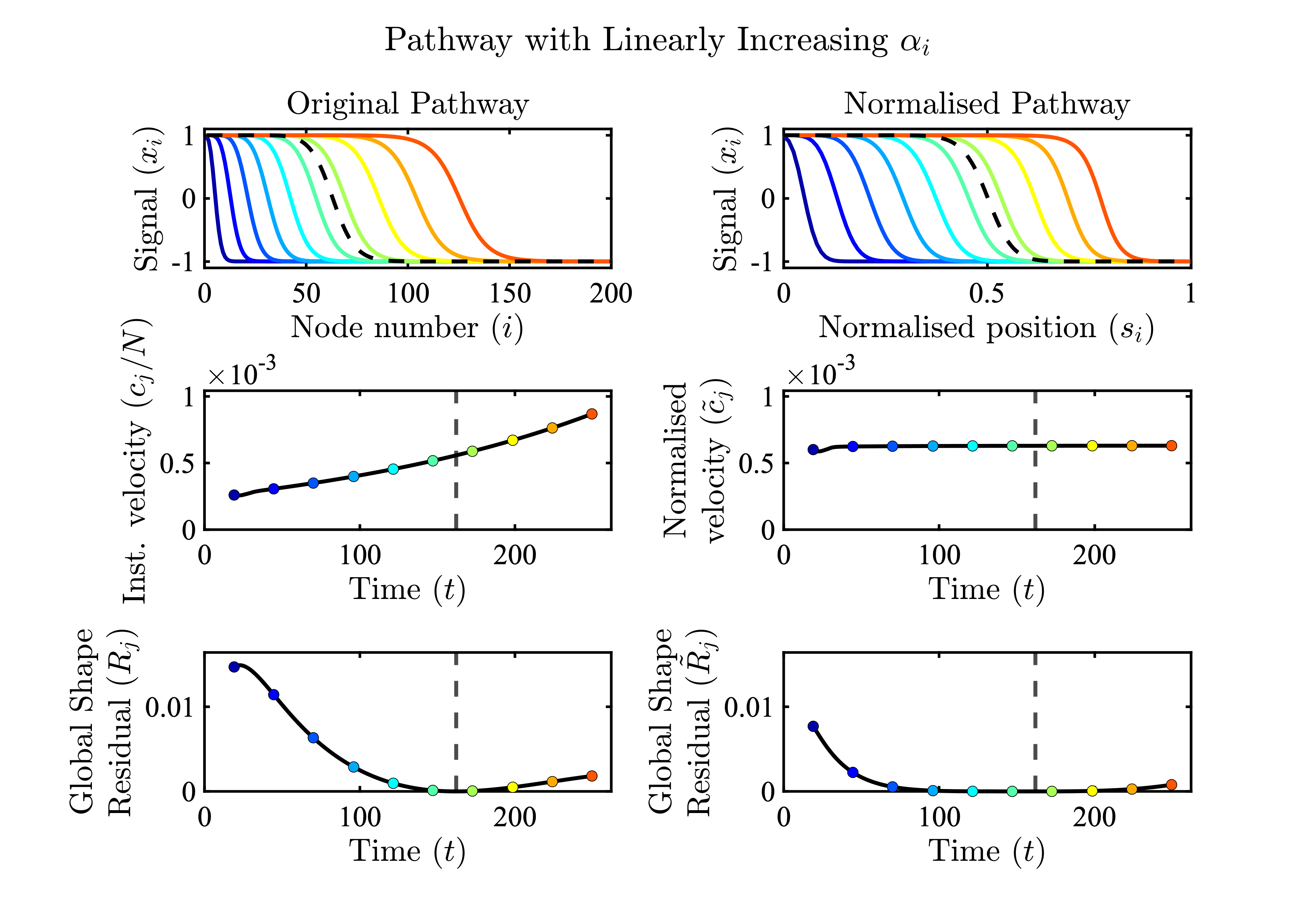} 
    \caption{Signal propagation dynamics with linearly increasing timescale parameter $\alpha_i \in [1, 5]$. Base parameters are $B = 100$, $\phi = 0$, $N = 200$, under boundary input $x_0 = 1$. The color gradient of the markers corresponds to time evolution, from early (dark blue) to late (red). The dashed vertical line marks the reference time $t_J$ at which the global shape residual is zero. \textit{Left:} The original spatial domain shows uneven spacing and a monotonically increasing measured velocity $c_j$ (Eq.~\ref{eq:cj}), reflecting continuous acceleration. \textit{Right:} The rescaled coordinate $s_i$ (Eq.~\ref{eq:rescaling}) produces evenly spaced wavefronts and a remarkably stable velocity $\tilde{c}_j$ (Eq.~\ref{eq:cjtilde}), confirming intrinsic wave dynamics remain uniform once the spatially varying timescale is absorbed.}
    \label{fig:active_alpha}
\end{figure}

We assess the effectiveness and limits of the heuristic assertion that each edge approximately confers its own independent delay to the pathway signal (that in the normalized pathway approximately uniform traveling waves are observed) by first examining pathways with monotonically varying parameters. Figure~\ref{fig:active_alpha} considers a pathway in which the timescale parameter $\alpha_i$ increases linearly with edge index $i$, while operating at a high saturation threshold ($B = 100$) and under an unbiased chemical mechanism ($\phi = 0$). Following our established framework, we apply the steady-state boundary input $x_0 = 1$ to a pathway initially at rest so that the stationary distribution is uniform. Because successive nodes respond progressively faster to their upstream input, the wavefront in the original domain undergoes continuous acceleration and wavefront dilation. This is shown by the colored markers progressing from early times (dark blue) to late times (red). We observe early transient discrepancies in the shape residual as the initial condition forms into a traveling wave, eventually converging toward zero at the reference time $t_J \approx 160$ (indicated by the vertical dashed line). 

However, when plotted against the rescaled coordinate $s_i$, the wavefronts become evenly spaced and the velocity profile stabilizes to a near-constant value ($\tilde{c}_j$) after the transient transition from the initial condition. In this high $B$ regime, the system kinetics approach the mass-action limit, resulting in inherently linear rates of transition across the wavefront as a result the spatial rescaling cleanly factors out the parameter-driven acceleration to recover the uniform intrinsic wave dynamics and remarkably fixing the shape heterogeneities. However, as the wavefronts get wider at larger values of $\alpha$ high order effects start to creep into the the distortion of the of the wavefront which can be seen in a small late stage increase in Global Shape Residual in the normalized pathway and this high order effect becomes more important for smaller $B$ values which is seen in the proceeding figures.

\begin{figure}[h]
    \centering
    \includegraphics[width=1\textwidth]{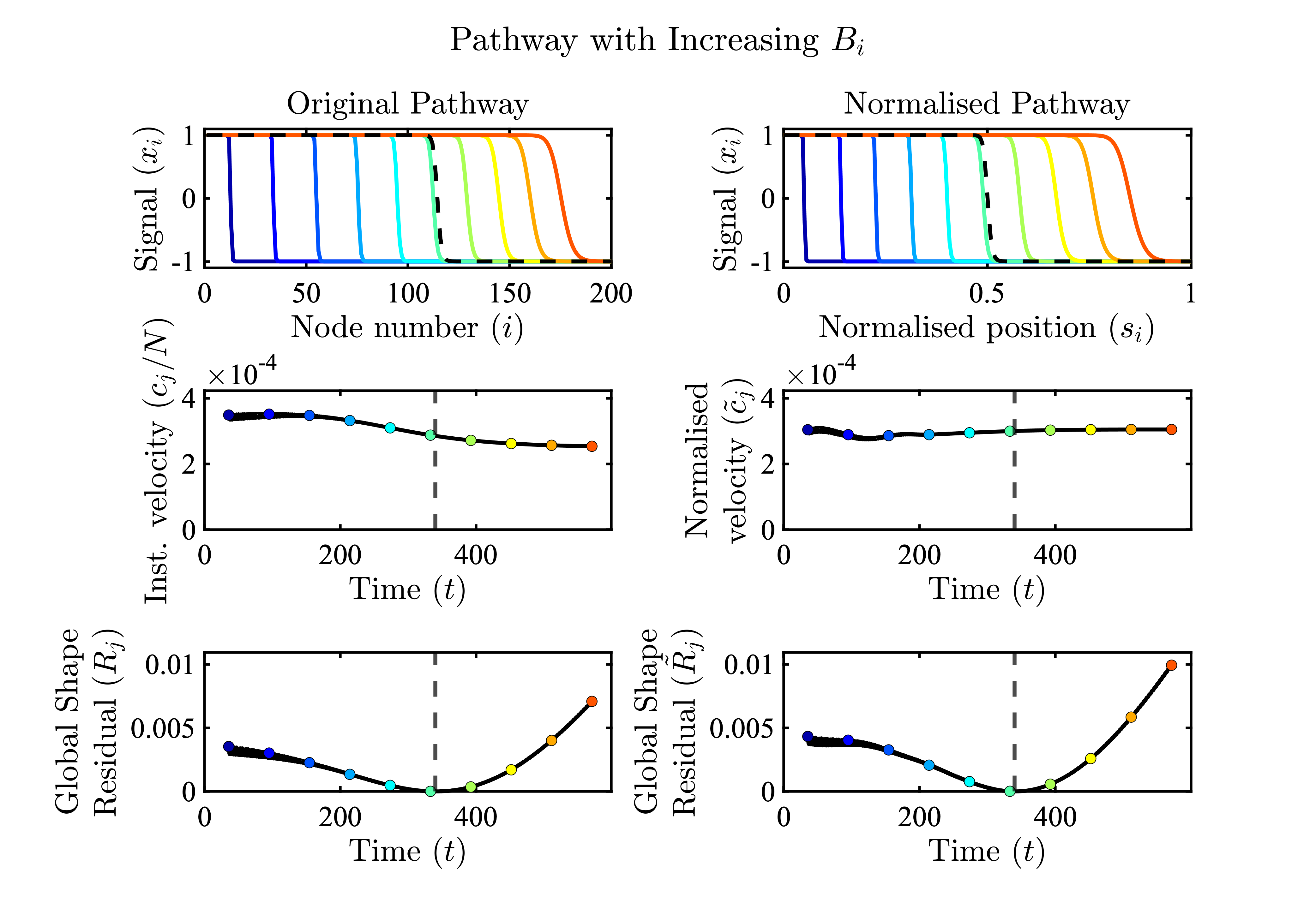} 
    \caption{Signal propagation dynamics with logarithmically increasing saturation threshold $B_i \in [1.02, 201]$. Base parameters are $\alpha = 1$, $\phi = 0$, $N = 200$. \textit{Left:} The velocity $c_j$ remains relatively flat compared to timescale variations, though the shape residual displays transient boundary effects early on (blue markers). \textit{Right:} Normalization successfully flattens the velocity $\tilde{c}_j$, but stretching the inherently wide wavefronts associated with low $B$ domains exacerbates the global shape residual $\tilde{R}_j$ after the reference time (dashed line, $t_J \approx 330$).}
    \label{fig:active_B}
\end{figure}

In contrast, Figure~\ref{fig:active_B} illustrates the effects of increasing the nonlinearity parameter $B_i$ logarithmically from $1 + 2\times10^{-2}$ (immediate saturation) to $1+2\times10^2$ (extremely high saturation). Because $B_i$ dictates the severity of the nonlinearity, variations in this parameter have comparatively little effect on the real pathway waveform's propagation speed (yielding a relatively flat $c_j$ compared to the $\alpha$ variation), but as $B_i$ increases the speed decreases slightly while the wavefront dilates. Consequently, when the domain is stretched via Equation~\ref{eq:rescaling} to keep the velocity approximately steady, it artificially dilates the wavefront even more leading to a global shape residual ($\tilde{R}_j$) that is slightly \textit{worse} than the original domain, highlighting a core limitation of the heuristic framework when applied to wide, non-mass-action signals. It should be noted that this distortion is small despite 4 orders of magnitude in variation of $B_i$.

\begin{figure}[h]
    \centering
    \includegraphics[width=1\textwidth]{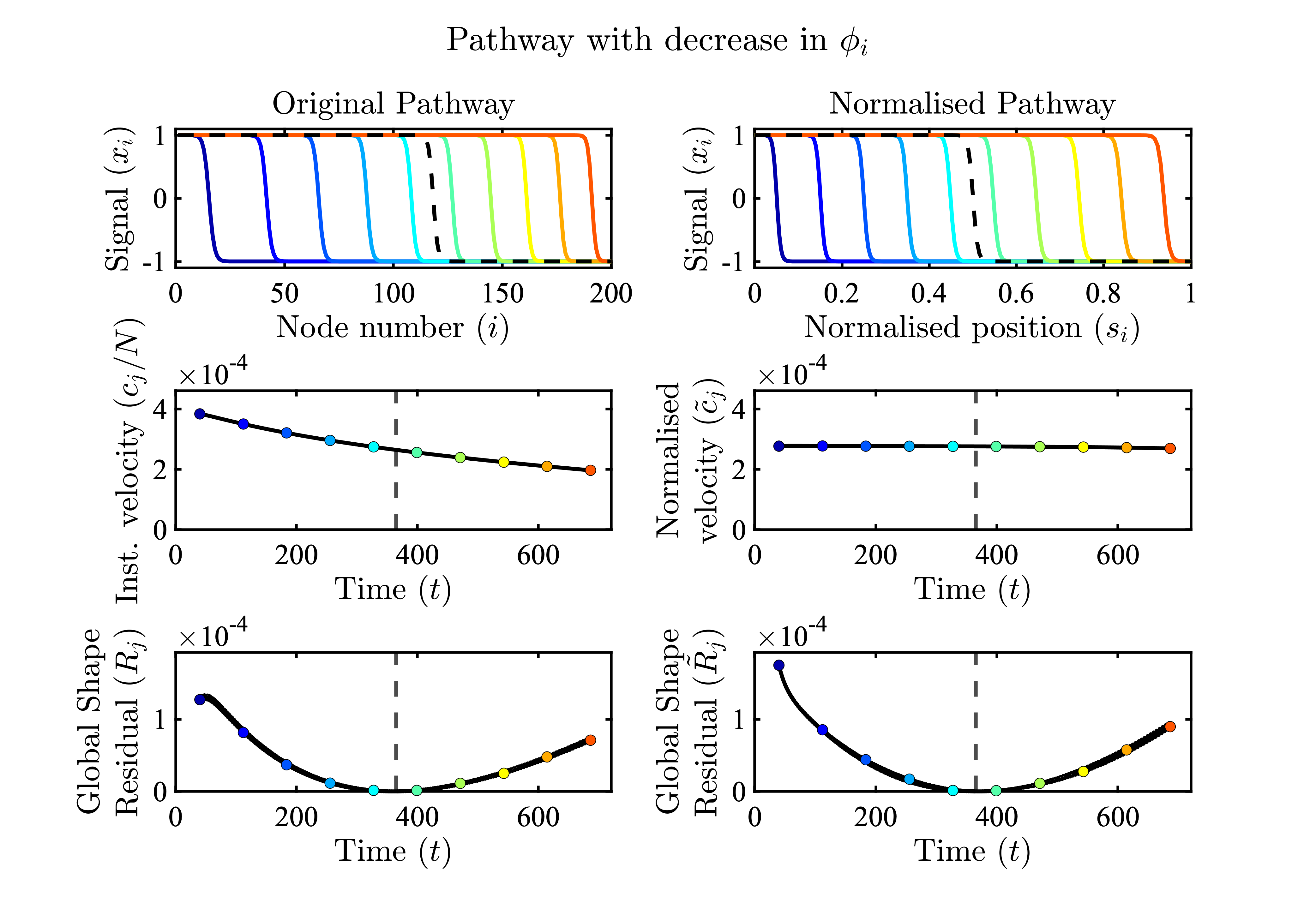} 
    \caption{Signal propagation dynamics with linearly decreasing asymmetry parameter $\phi_i$. Base parameters are $\alpha = 1$, $B = 5$, $N = 200$. Both original and normalized residuals ($R_j$, $\tilde{R}_j$) capture early transient effects (dark blue markers). The requirement of a low $B_0$ to permit meaningful $\phi_i$ dynamics results in a wide wavefront that undergoes structural distortion upon normalization. High-frequency noise at late times (red markers) indicates numerical instabilities as the leading edge of the wide wavefront interacts with the terminal boundary before the reference time $t_J \approx 375$ (dashed line) is fully resolved.}
    \label{fig:active_phi}
\end{figure}

This distortion effect is similarly present when varying the asymmetry parameter $\phi_i$ (Figure~\ref{fig:active_phi}). Here, the pathway is initialized with a moderately low saturation limit ($B_0 = 5$) to allow the decreasing $\phi_i$ to meaningfully impact the local kinetics. Again, the wide wavefronts (in this case early in time) suffer geometric distortion upon spatial normalization. It is important to note, however that monotonically varying each of the three parameters in a controlled sense and implementing the normalized pathway gives rise to a very stable wavespeed which gives robust predictability around the lag time between input and downstream changes. To push this further, we will now investigate how robust the travelling wave is under normalization of unstructured stochastic edge parameters. 

\subsection{Pathways with random parameter variation}

To evaluate the robustness of the rescaling framework (the normalised pathway) in biologically realistic environments, we introduce stochastic heterogeneity. In real biochemical cascades,  kinetic parameters fluctuate randomly from edge to edge. To model this, we assign parameters drawn from heavy-tailed lognormal distributions, ensuring all rates remain strictly positive. We set the base activation rates $\alpha_i \sim \mathrm{Lognormal}(\ln \alpha_0, \sigma)$, where $\alpha_0 = 1$ is the baseline median. To guarantee the nonlinearity parameter $\beta_i$ strictly exceeds $1$, we apply the lognormal distribution to the shifted variable $X_i = \beta_i - 1$, such that $X_i \sim \mathrm{Lognormal}(\ln(\beta_0 - 1), \sigma)$ with a baseline median $\beta_0 = 5$. The parameter $\sigma$ acts as the singular tuning dial for the degree of spatial heterogeneity across the pathway.

\begin{figure}[htbp]
    \centering
    \includegraphics[width=1\textwidth]{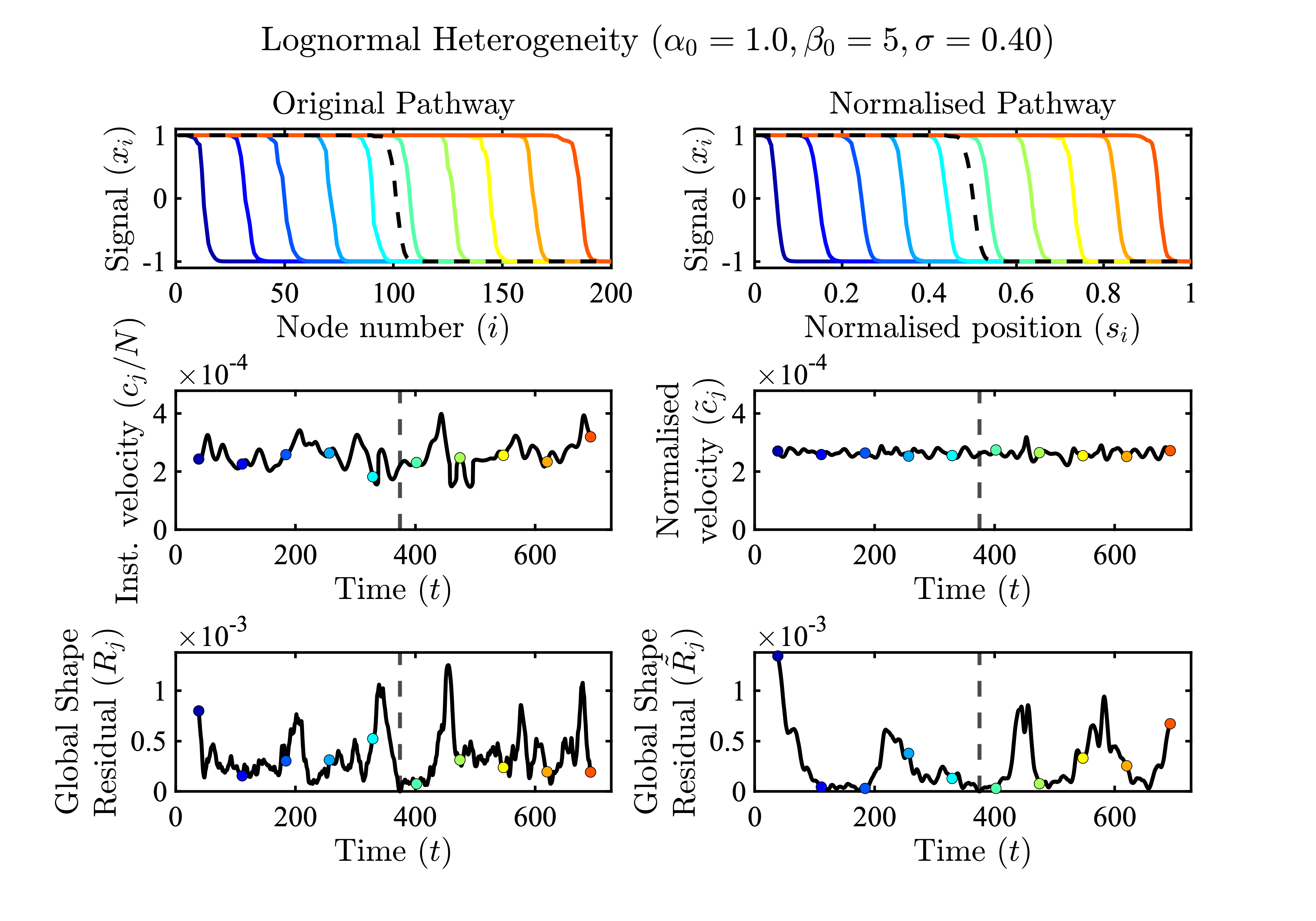}
    \caption{Signal propagation through a stochastic pathway ($\alpha_0=1, \beta_0=5, \sigma=0.4, \phi_i=0$). The initial transient phase is omitted, with velocity and residual metrics recorded only after the wave passes $s=0.15$. Random parameter variations induce heavy fluctuations in the original instantaneous velocity (left column). Rescaling the spatial domain (right column) absorbs these kinetic variations, resulting in a highly uniform normalised velocity $\tilde{c}_j$ while maintaining low global shape residuals, successfully recovering a predictable travelling wave structure.}
    \label{fig:stoch_single}
\end{figure}

Figure~\ref{fig:stoch_single} illustrates the framework applied to a single stochastic realisation with moderate variability ($\sigma = 0.4$, $\phi_i = 0$). To evaluate the steady-state travelling wave properties accurately, we omit the initial transient phase by restricting our metric calculations to the temporal window after the wave has traversed the first 15\% of the normalised domain ($s \ge 0.15$). In the original spatial coordinate system, random parameter anomalies act as local bottlenecks or accelerants, causing the wave to stutter and the instantaneous velocity $c_j/N$ to heavily fluctuate. By contrast, the reciprocal-velocity rescaling successfully absorbs these irregularities. The normalised velocity profile $\tilde{c}_j$ is substantially flattened, demonstrating that the transformation successfully decouples the intrinsic wave dynamics from parameter-driven noise. Furthermore, these gains in velocity uniformity do not come at the cost of the wave's structural integrity; the global shape residual $\tilde{R}_j$ of the normalised pathway remains comparable to, or slightly lower than, the original $R_j$, indicating a smoothly translating wavefront.

\begin{figure}[htbp]
    \centering
    \includegraphics[width=1\textwidth]{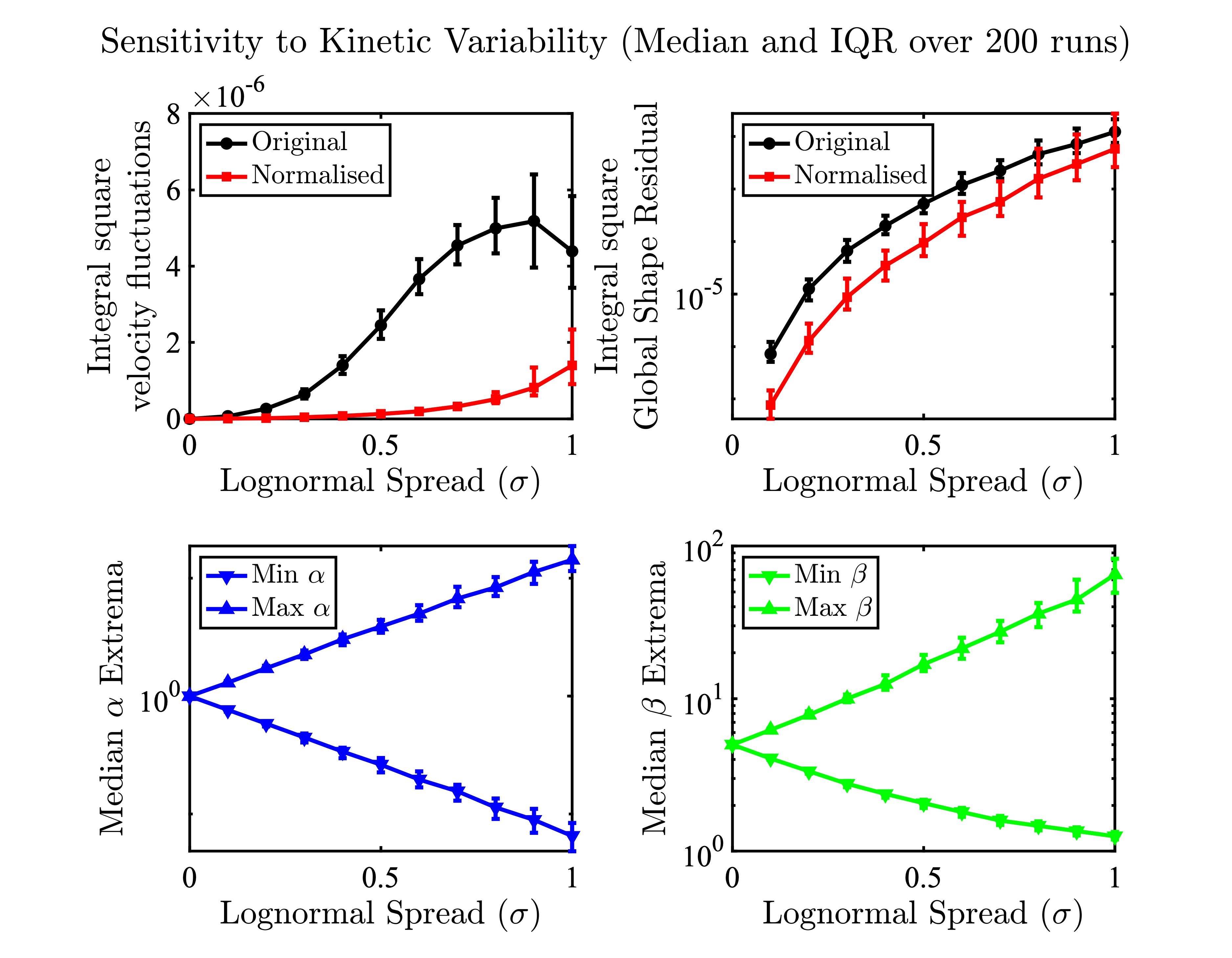}
    \caption{Statistical sensitivity of wave propagation to kinetic parameter variability ($\sigma$), evaluated over 200 realisations per point. Metrics are integrated over the post-transient time window. (a) Median (and inter-quartile range) Velocity Integral Square Error (VISE) for both original and normalised pathways --a measure of velocity fluctuation. Normalisation consistently suppresses velocity fluctuations at low to moderate $\sigma$. The artificial drop in original VISE at high $\sigma$ is due to pathway fragmentation. (b) Median (and inter-quartile range) Residual Integral Square Error (RISE) for both original and normalised pathways --a measure of wavefront distortion. Normalisation maintains or slightly improves wavefront shape stability, though variance increases at extreme $\sigma$ as continuous fitting fails across severe discrete bottlenecks. (c, d) The average extrema of $\alpha_i$ and $\beta_i$ respectively, demonstrating the immense parameter spread (exceeding two orders of magnitude) driving the pathway fragmentation at $\sigma = 1$.}
    \label{fig:stoch_sweep}
\end{figure}

To systematically quantify the efficacy of the framework across different noise regimes, we condense the temporal metrics into a single summary statistic: the Integral Square Error (ISE). Over the valid post-transient time window $[t_{\mathrm{start}}, t_{\mathrm{end}}]$, we define the Velocity ISE (VISE) as the integral of the squared deviations of the velocity from its time-averaged mean, strongly penalising large velocity fluctuations. Similarly, the Residual ISE (RISE) is calculated by integrating the squared global shape residual over time. We reinitialise 200 distinct random pathways at sequential levels of $\sigma \in [0, 1]$, extracting the median and interquartile range (IQR) of these integral metrics (Figure~\ref{fig:stoch_sweep}).

The statistical sweep reveals several key insights into the operating limits of the rescaling methodology. Firstly, Figure~\ref{fig:stoch_sweep}a confirms that for low to moderate heterogeneity ($\sigma \lesssim 0.6$), the normalised pathway dramatically reduces velocity fluctuations compared to the original domain. Secondly, the shape residual metric (Figure~\ref{fig:stoch_sweep}b) demonstrates that the normalisation procedure slightly improves or preserves wavefront stability in this regime.

However, the behaviour at high heterogeneity ($\sigma > 0.6$) is counter-intuitive: the original pathway's velocity fluctuations reach a maximum and then decrease, despite the pathway heterogeneity becoming maximally noisy. This decrease in the original VISE does not indicate a return to smooth wave propagation. Rather, it is a symptom of extreme parameter divergence. As shown in Figures~\ref{fig:stoch_sweep}c and~\ref{fig:stoch_sweep}d, at $\sigma = 1$, the minimum and maximum kinetic parameters span over two orders of magnitude. 

In this extreme regime, the pathway contains a sufficient number of edges with vastly elevated speeds (high $\alpha$), interspersed with edges possessing severely stunted speeds. The ultra-fast edges propagate information so rapidly that clusters of these nodes essentially collapse into a single state, acting as a single node in quasi-equilibrium. Conversely, the ultra-slow edges dictate the rate of the entire system. Consequently, the pathway undergoes a functional fragmentation. The effectively shortened pathway is dominated entirely by the slow bottlenecks, filtering out the fast dynamics and artificially reducing the velocity variance over time. The normalised pathway, which attempts to stretch the coordinate space continuously across these severe discrete bottlenecks, exhibits rising shape and velocity variance, identifying the functional limit of assuming continuous wave-like propagation in highly discrete, fragmented cascades.

\section{Discussion and conclusion}

In this work, we have developed a comprehensive mathematical framework for characterising signal propagation in non-linear, heterogeneous biochemical pathways. By modelling these pathways as cascades of discrete interactions governed by Michaelis-Menten-type kinetics, we established exact steady-state distributions and identified conditions supporting travelling wave solutions. Crucially, we demonstrated that spatial heterogeneities whether smoothly graded or highly stochastic—distort classical wave propagation, leading to velocity stalling and wavefront deformation. To mitigate this, we introduced a heuristic reciprocal-velocity spatial rescaling. This coordinate transformation successfully absorbs parameter-driven structural noise, smoothing wave velocities and preserving structural wave profiles without requiring bespoke parameter tuning or continuous limits.

The sensitivity sweep in stochastic environments defines the clear operating envelope of our framework. The reciprocal-velocity normalisation is profoundly effective in the moderate heterogeneity regime ($\sigma \in [0.1, 0.5]$), which aligns with biologically plausible levels of intrinsic noise in gene expression and enzyme concentrations. However, when the parameter distribution becomes extremely heavy-tailed ($\sigma \to 1$), the framework’s assumptions break down. In this regime, the system is dominated by severe bottlenecks, stretching the reciprocal-velocity weights excessively and preventing a continuous travelling wave shape from fitting the discrete, disjointed signal transitions. 

Interestingly, this breakdown offers a profound insight into how mathematical modellers can abstract and reduce large-scale biological networks in a systematic way. Formulating a model that explicitly describes every sequential node in a large pathway (or any serialised biological process) often leads to analytical intractability, computational stiffness, and overfitting to sparse experimental data. Our analysis theoretically justifies a principled model reduction strategy based on kinetic thresholds:
\begin{enumerate}
    \item \textbf{Ultra-fast edges ($\alpha_i \gg 1$):} These interactions rapidly reach equilibrium. As we observed in the high-$\sigma$ regime, these clusters of nodes effectively coarsen, justifying their algebraic collapse into single functional nodes within a reduced model.
    \item \textbf{Moderate wave-like segments:} Contiguous sequences of nodes lacking severe bottlenecks exhibit predictable, smoothed travelling wave behaviour under our normalisation. These segments can be abstracted mathematically by replacing explicit sequential ODEs with single Delay Differential Equations (DDEs), where the fixed signal transit time is calculated explicitly using our theoretical asymptotic speed heuristic.
    \item \textbf{Ultra-slow edges ($\alpha_i \ll 1$):} These segments act as strict bottlenecks. Modellers must break the wave abstraction at these locations, retaining the slow transitions as explicitly modelled discrete ODE species.
\end{enumerate}
By separating a pathway into these three dynamical classes, modellers can construct hybrid discrete-delay models that drastically reduce system dimensionality while rigorously preserving the precise temporal dynamics and delays inherent to signal transduction.

\bibliography{sn-bibliography}

@article{fisher1937wave,
  title={The wave of advance of advantageous genes},
  author={Fisher, Ronald Aylmer},
  journal={Annals of eugenics},
  volume={7},
  number={4},
  pages={355--369},
  year={1937},
  publisher={Wiley Online Library}
}

@article{McKeithan1995KineticProofreading,
  author    = {McKeithan, Thomas W.},
  title     = {Kinetic proofreading in T-cell receptor signal transduction},
  journal   = {Proceedings of the National Academy of Sciences of the United States of America},
  year      = {1995},
  volume    = {92},
  number    = {11},
  pages     = {5042--5046},
  doi       = {10.1073/pnas.92.11.5042}
}

@article{hurtado2019generalizations,
  title={Generalizations of the ‘Linear Chain Trick’: incorporating more flexible dwell time distributions into mean field ODE models},
  author={Hurtado, Paul J and Kirosingh, Adam S},
  journal={Journal of mathematical biology},
  volume={79},
  number={5},
  pages={1831--1883},
  year={2019},
  publisher={Springer}
}

@article{scarabel2024equations,
  title={Equations with infinite delay: pseudospectral discretization for numerical stability and bifurcation in an abstract framework},
  author={Scarabel, Francesca and Vermiglio, Rossana},
  journal={SIAM Journal on Numerical Analysis},
  volume={62},
  number={4},
  pages={1736--1758},
  year={2024},
  publisher={SIAM}
}

@book{Murray2002MathematicalBiologyI,
  author    = {Murray, James D.},
  title     = {Mathematical Biology. I: An Introduction},
  edition   = {3},
  year      = {2002},
  publisher = {Springer},
  address   = {New York},
  series    = {Interdisciplinary Applied Mathematics},
  volume    = {17},
  isbn      = {978-0387952239}
}

@article{Glass2021NonlinearDelayMotifs,
  author    = {Glass, Daniel S. and Jang, Sung H. and Fraser, Hunter B.},
  title     = {Nonlinear delay differential equations and their application to modeling biological network motifs},
  journal   = {Nature Communications},
  year      = {2021},
  volume    = {12},
  number    = {1},
  pages     = {1788},
  doi       = {10.1038/s41467-021-21700-8},
  url       = {https://doi.org/10.1038/s41467-021-21700-8}
}

@article{Schwen2021NonlinearSignaling,
  author       = {Schwen, Peter and Geiger, Tam{\'a}s},
  title        = {Modeling the Nonlinear Dynamics of Intracellular Signaling Networks},
  journal      = {Bio-protocol},
  year         = {2021},
  volume       = {11},
  number       = {14},
  pages        = {e4089},
  doi          = {10.21769/BioProtoc.4089},
  url          = {https://bio-protocol.org/e4089},
  note         = {Protocol on nonlinear ODE and PDE models of intracellular signaling dynamics, including multistability and oscillations},
}

@incollection{Othmer1997SignalTransduction,
  author       = {Othmer, Hans G.},
  title        = {Signal Transduction and Second Messenger Systems},
  booktitle    = {Case Studies in Mathematical Modeling: Ecology, Physiology, and Cell Biology},
  editor       = {Othmer, Hans G. and Adler, Frederick R. and Lewis, Mark A. and Dallon, John C.},
  publisher    = {Prentice Hall},
  address      = {Englewood Cliffs, NJ},
  year         = {1997},
  pages        = {101--128},
  isbn         = {0-13-574039-8},
  note         = {Classic mathematical treatment of signal transduction dynamics with multiple qualitative behaviors},
}

@article{kolmogorov1982study,
  title={Study of the diffusion equation with growth of the quantity of matter and its applications to biological problems},
  author={Kolmogorov, A},
  journal={Applicable Mathematics to Non-physical Phenomena},
  year={1982},
  publisher={John Wiley Sons}
}

@article{hunter2000signaling,
  title={Signaling—2000 and beyond},
  author={Hunter, Tony},
  journal={Cell},
  volume={100},
  number={1},
  pages={113--127},
  year={2000},
  publisher={Elsevier}
}

@article{pearson2001mitogen,
  title={Mitogen-activated protein (MAP) kinase pathways: regulation and physiological functions},
  author={Pearson, Gray and Robinson, Fred and Beers Gibson, Tara and Xu, Bing-e and Karandikar, Mahesh and Berman, Kevin and Cobb, Melanie H},
  journal={Endocrine reviews},
  volume={22},
  number={2},
  pages={153--183},
  year={2001},
  publisher={Oxford University Press}
}

@article{clevers2012wnt,
  title={Wnt/$\beta$-catenin signaling and disease},
  author={Clevers, Hans and Nusse, Roel},
  journal={Cell},
  volume={149},
  number={6},
  pages={1192--1205},
  year={2012},
  publisher={Elsevier}
}

@article{huang1996ultrasensitivity,
  title={Ultrasensitivity in the mitogen-activated protein kinase cascade.},
  author={Huang, Chi-Ying and Ferrell Jr, James E},
  journal={Proceedings of the National Academy of Sciences},
  volume={93},
  number={19},
  pages={10078--10083},
  year={1996},
  publisher={National Acad Sciences}
}

@article{kholodenko2000negative,
  title={Negative feedback and ultrasensitivity can bring about oscillations in the mitogen-activated protein kinase cascades},
  author={Kholodenko, Boris N},
  journal={European journal of biochemistry},
  volume={267},
  number={6},
  pages={1583--1588},
  year={2000},
  publisher={Wiley Online Library}
}

@article{jayathilaka2024two,
  title={Two wrongs do not make a right: the assumption that an inhibitor acts as an inverse activator},
  author={Jayathilaka, Chathranee and Araujo, Robyn and Nguyen, Lan and Flegg, Mark},
  journal={Journal of Mathematical Biology},
  volume={89},
  number={2},
  pages={26},
  year={2024},
  publisher={Springer}
}

@article{heinrich2002mathematical,
  title={Mathematical models of protein kinase signal transduction},
  author={Heinrich, Reinhart and Neel, Benjamin G and Rapoport, Tom A},
  journal={Molecular cell},
  volume={9},
  number={5},
  pages={957--970},
  year={2002},
  publisher={Elsevier}
}

@article{martiel1987model,
  title={A model based on receptor desensitization for cyclic AMP signaling in Dictyostelium cells},
  author={Martiel, Jean-Louis and Goldbeter, Albert},
  journal={Biophysical journal},
  volume={52},
  number={5},
  pages={807--828},
  year={1987},
  publisher={Elsevier}
}

@incollection{aronson2006nonlinear,
  title={Nonlinear diffusion in population genetics, combustion, and nerve pulse propagation},
  author={Aronson, Donald G and Weinberger, Hans F},
  booktitle={Partial Differential Equations and Related Topics: Ford Foundation Sponsored Program at Tulane University, January to May, 1974},
  pages={5--49},
  year={2006},
  publisher={Springer}
}

@article{Rios2021,
author = "Karina Islas Rios",
title = "{Systematic characterisation of network structures underlying adaptive resistance and susceptibility to sequential combination treatment in cancer}",
year = "2021",
month = "12",
url = "https://bridges.monash.edu/articles/thesis/Systematic_characterisation_of_network_structures_underlying_adaptive_resistance_and_susceptibility_to_sequential_combination_treatment_in_cancer/17209520",
doi = "10.26180/17209520.v1"
}

@article{o2011tunable,
  title={Tunable signal processing in synthetic MAP kinase cascades},
  author={O'Shaughnessy, Ellen C and Palani, Santhosh and Collins, James J and Sarkar, Casim A},
  journal={Cell},
  volume={144},
  number={1},
  pages={119--131},
  year={2011},
  publisher={Elsevier}
}

\begin{appendices}


\end{appendices}

\end{document}